# Joint continuity of the local times of fractional Brownian sheets

Antoine Ayache[a], Dongsheng Wu[b,1] and Yimin Xiao[c,2]

[a]*U.M.R. CNRS 8524, Laboratoire Paul Painleve, Bat. M2, Université Lille 1, 59655 Villenuve d'ascq Cedex, France and U.M.R. CNRS 8179, LEM, IAE de Lille, 104, Avenue du peuple Belge 59043 Lille Cedex, France.*
*E-mail: Antoine.Ayache@math.univ-lille1.fr*

[b]*Department of Mathematical Sciences, 339A Madison Hall, The University of Alabama in Huntsville, Huntsville, AL 35899, USA. E-mail: wudongsh@msu.edu*

[c]*Department of Statistics and Probability, A-413 Wells Hall, Michigan State University, East Lansing, MI 48824, USA.*
*E-mail: xiao@stt.msu.edu*

**Abstract.** Let $B^H = \{B^H(t), t \in \mathbb{R}_+^N\}$ be an $(N,d)$-fractional Brownian sheet with index $H = (H_1, \ldots, H_N) \in (0,1)^N$ defined by $B^H(t) = (B_1^H(t), \ldots, B_d^H(t))$ ($t \in \mathbb{R}_+^N$), where $B_1^H, \ldots, B_d^H$ are independent copies of a real-valued fractional Brownian sheet $B_0^H$. We prove that if $d < \sum_{\ell=1}^N H_\ell^{-1}$, then the local times of $B^H$ are jointly continuous. This verifies a conjecture of Xiao and Zhang (*Probab. Theory Related Fields* **124** (2002)).

We also establish sharp local and global Hölder conditions for the local times of $B^H$. These results are applied to study analytic and geometric properties of the sample paths of $B^H$.

**Résumé.** Désignons par $B^H = \{B^H(t), t \in \mathbb{R}_+^N\}$ le $(N,d)$-drap Brownien fractionnaire de paramètre $H = (H_1, \ldots, H_N) \in (0,1)^N$ défini par $B^H(t) = (B_1^H(t), \ldots, B_d^H(t))$ ($t \in \mathbb{R}_+^N$), où $B_1^H, \ldots, B_d^H$ sont des copies indépendantes du drap Brownien fractionnaire à valeurs réelles $B_0^H$. Nous montrons que le temps local de $B^H$ est bicontinu lorsque $d < \sum_{\ell=1}^N H_\ell^{-1}$. Cela résout une conjecture de Xiao et Zhang (*Probab. Theory Related Fields* **124** (2002)). Nous obtenons aussi des résultats fins concernant la régularité Hölderienne, locale et globale, du temps local. Ces résultats nous permettent d'étudier certaines propriétés analytiques et géométriques des trajectoires de $B^H$.

[1]Research partially supported by NSF Grant DMS-0417869.
[2]Research partially supported by NSF Grant DMS-0404729.

## 1. Introduction

For a given vector $H = (H_1, \ldots, H_N) \in (0,1)^N$, a real-valued fractional Brownian sheet $B_0^H = \{B_0^H(t), t \in \mathbb{R}_+^N\}$ with index $H$ is a centered Gaussian random field with covariance function given by

$$\mathbb{E}[B_0^H(s)B_0^H(t)] = \prod_{\ell=1}^{N} \frac{1}{2}(s_\ell^{2H_\ell} + t_\ell^{2H_\ell} - |s_\ell - t_\ell|^{2H_\ell}), \quad s,t \in \mathbb{R}_+^N. \tag{1.1}$$

It follows from (1.1) that $B_0^H(t) = 0$ a.s. for every $t \in \partial \mathbb{R}_+^N$, where $\partial \mathbb{R}_+^N$ denotes the boundary of $\mathbb{R}_+^N$.

We will make use of the following stochastic integral representation of $B_0^H$ (cf. [2]):

$$B_0^H(t) = \kappa_H^{-1} \int_{-\infty}^{t_1} \cdots \int_{-\infty}^{t_N} \prod_{\ell=1}^{N} g_{H_\ell}(t_\ell, s_\ell) W(\mathrm{d}s), \tag{1.2}$$

where $W = \{W(s), s \in \mathbb{R}^N\}$ is a standard real-valued Brownian sheet and where, for every $\ell = 1, \ldots, N$,

$$g_{H_\ell}(t_\ell, s_\ell) = ((t_\ell - s_\ell)_+)^{H_\ell - 1/2} - ((-s_\ell)_+)^{H_\ell - 1/2}.$$

In the above, $a_+ = \max\{a, 0\}$ for all $a \in \mathbb{R}$ and $\kappa_H$ is the normalization constant given by

$$\kappa_H^2 = \int_{-\infty}^{1} \cdots \int_{-\infty}^{1} \left[\prod_{\ell=1}^{N} g_{H_\ell}(1, s_\ell)\right]^2 \mathrm{d}s.$$

Note that if $H_{\ell_0} = 1/2$ for some $\ell_0$, then we assume that $g_{H_{\ell_0}}(t_{\ell_0}, s_{\ell_0}) = \mathbb{1}_{[0, t_{\ell_0}]}(s_{\ell_0})$, where $\mathbb{1}_{[0, t_{\ell_0}]}$ is the indicator of the interval $[0, t_{\ell_0}]$.

Let $B_1^H, \ldots, B_d^H$ be $d$ independent copies of $B_0^H$. Then the Gaussian random field $B^H = \{B^H(t): t \in \mathbb{R}_+^N\}$ with values in $\mathbb{R}^d$ defined by

$$B^H(t) = (B_1^H(t), \ldots, B_d^H(t)), \quad \forall t \in \mathbb{R}_+^N, \tag{1.3}$$

is called an $(N, d)$-fractional Brownian sheet with index $H = (H_1, \ldots, H_N)$.

Note that if $N = 1$, then $B^H$ is a fractional Brownian motion in $\mathbb{R}^d$ with Hurst index $H_1 \in (0, 1)$; if $N > 1$ and $H_1 = \cdots = H_N = 1/2$, then $B^H$ is the $(N, d)$-Brownian sheet. However, when $H_1, \ldots, H_N$ are not the same, $B^H$ is anisotropic and has the following operator-self-similarity (this can be verified easily using (1.1)): For any $N \times N$ diagonal matrix $A = (a_{ij})$ with $a_{ii} = a_i > 0$ for all $1 \leq i \leq N$ and $a_{ij} = 0$ if $i \neq j$, we have

$$\{B^H(At), t \in \mathbb{R}_+^N\} \stackrel{d}{=} \left\{\prod_{j=1}^{N} a_j^{H_j} B^H(t), t \in \mathbb{R}_+^N\right\}, \tag{1.4}$$

where $X \stackrel{d}{=} Y$ means that the two processes have the same finite dimensional distributions. These features of $B^H$ make it a possible model for bone structure [8] and aquifer structure in hydrology [4].

Many authors have studied various properties of fractional Brownian sheets. See, for example, [3, 11, 21, 24, 29, 33] and the references therein for further information. This paper is concerned with regularity of the local times of an $(N, d)$-fractional Brownian sheet $B^H$. After having proved that a necessary and sufficient condition for the existence of $L^2(\mathbb{P} \times \lambda_d)$ local times of $B^H$ is $d < \sum_{\ell=1}^{N} \frac{1}{H_\ell}$, Xiao and Zhang [33] give a sufficient condition for the joint continuity of the local times. However, their sufficient condition is not sharp and they have conjectured that $B^H$ has jointly continuous local times whenever the condition $d < \sum_{\ell=1}^{N} \frac{1}{H_\ell}$ is satisfied. The main objective of this paper is to verify this conjecture; see Theorem 3.1. The new ingredients for proving this result is the property of sectorial local nondeterminism of $B_0^H$ established in [29] (see Lemma 3.2) and a similar result for the fractional Liouville sheet proved in Section 2. The results and techniques developed in this paper are applicable to more general anisotropic Gaussian random fields with the property of sectorial local nondeterminism; see [32] for further development.

The rest of this paper is organized as follows. In Section 2, we prove some basic results on the fractional Liouville sheets that will be useful to our arguments. In Section 3, we prove that the sufficient condition for the existence of $L^2(\mathbb{P} \times \lambda_d)$ local times of $B^H$ in [33] actually implies the joint continuity of the local times. This verifies their conjecture in Remark 4.11. Section 4 is on the local and uniform Hölder conditions for the local times and their implications to sample path properties of $B^H$. In particular, we derive some results on the Hausdorff measure of the level sets and on the Chung-type law of the iterated logarithm for the sample function $B^H(t)$. The latter improves Theorem 3 of [3].

We end the Introduction with some notation. Throughout this paper, the underlying parameter space is $\mathbb{R}^N$ or $\mathbb{R}^N_+ = [0, \infty)^N$. A parameter $t \in \mathbb{R}^N$ is written as $t = (t_1, \ldots, t_N)$, or as $\langle c \rangle$, if $t_1 = \cdots = t_N = c$. For any $s, t \in \mathbb{R}^N$ such that $s_j < t_j$ ($j = 1, \ldots, N$), we define the closed interval (or rectangle) $[s,t] = \prod_{j=1}^N [s_j, t_j]$. We will let $\mathcal{A}$ denote the class of all closed intervals $I \subset (0, \infty)^N$. We always write $\lambda_m$ for Lebesgue's measure on $\mathbb{R}^m$, and use $\langle \cdot, \cdot \rangle$ and $|\cdot|$ to denote the ordinary scalar product and the Euclidean norm in $\mathbb{R}^m$ respectively, no matter the value of the integer $m$.

An unspecified positive and finite constant will be denoted by $c$, which may not be the same in each occurrence. More specific constants in Section $i$ are numbered as $c_{i,1}, c_{i,2}, \ldots$.

## 2. Fractional Liouville sheet

One of the main obstacles in studying the local times and other properties of fractional Brownian sheets is their complicated dependence structure. Unlike the Brownian sheet or fractional Brownian motion, fractional Brownian sheets have neither the property of independent increments nor the local nondeterminism.

To be more specific, we recall that fractional Brownian motion $Z^\alpha = \{Z^\alpha(t), t \in \mathbb{R}^N\}$ ($0 < \alpha < 1$) in $\mathbb{R}$ has the following property of strong local nondeterminism proved by Pitt [25]: For every interval $I \subseteq \mathbb{R}^N$, there exist positive constants $c_{2,1}$ and $r_0$ such that for all $t \in I$ and all $0 < r \leq \min\{|t|, r_0\}$,

$$\text{Var}(Z^\alpha(t) \mid Z^\alpha(s): s \in I, r \leq |s - t| \leq r_0) \geq c_{2,1} r^{2\alpha}. \tag{2.1}$$

This property has played important rôles in studying the local times and many other properties of $Z^\alpha$; see [31] and the references therein for more information. On the other hand, it is known that the Brownian sheet $W = \{W(t), t \in \mathbb{R}^N_+\}$ does not have the property of local nondeterminism. In order to see this, we consider the Brownian sheet with $N = 2$ and $I = [0,1]^2$. For any constant $\varepsilon \in (0,1)$, let $T \subseteq I$ be an interval with side-length $\varepsilon$. Let $t$ denote the upper-right vertex of $T$ and let $s^1, s^2, s^3$ be other vertices of $T$. For example, $t = (1,1)$, $s^1 = (1-\varepsilon, 1), s^2 = (1, 1-\varepsilon)$ and $s^3 = (1-\varepsilon, 1-\varepsilon)$. Then $|t - s^j| \geq \varepsilon$ for $j = 1, 2, 3$. Considering the increment of $W$ over the square $T$, we see that $\text{Var}(W(t)|W(s^1), W(s^2), W(s^3)) \leq \varepsilon^2$. Hence the Brownian sheet $W$ does not satisfy (2.1) (this also proves that the fractional Brownian sheet $B^H_0$ is not locally nondeterministic). This is the main reason why, in most literature, the methods for studying various properties of the Brownian sheet are different from those for fractional Brownian motion. The property of independent increments of $W$ has been crucial in studying the local times and self-intersection local times of $W$; see [12, 27] and [22]. In solving an open problem in [22], Khoshnevisan and Xiao [19] showed that $W$ satisfies a type of sectorial local nondeterminism and applied this property to study geometric properties of the Brownian sheet by using methods that are reminiscent to those for fractional Brownian motion; see [18] for further applications of the sectorial local nondeterminism. Recently, Wu and Xiao [29] have extended several results in [18, 19] to fractional Brownian sheets.

In this paper we continue the above line of research and study the regularity of the local times of fractional Brownian sheets. To overcome the difficulty due to the lack of local nondeterminism of $B^H$, we will not only make use of the sectorial local nondeterministic property of $B^H$ established in [29] (see Lemma 3.2), but also the analogous properties of the so-called fractional Liouville sheet.

Given any vector $\alpha = (\alpha_1, \ldots, \alpha_N) \in (0, \infty)^N$, the centered Gaussian random field $X_0^\alpha = \{X_0^\alpha(t), t \in \mathbb{R}_+^N\}$ defined by

$$X_0^\alpha(t) = \int_{[0,t]} \prod_{\ell=1}^N (t_\ell - s_\ell)^{\alpha_\ell - 1/2} W(\mathrm{d}s), \quad t \in \mathbb{R}_+^N, \tag{2.2}$$

is called a fractional Liouville sheet with parameter $\alpha$. It is easy to see that, when $\alpha_1, \ldots, \alpha_N$ are not the same, $X_0^\alpha = \{X_0^\alpha(t) \colon t \in \mathbb{R}_+^N\}$ is an anisotropic Gaussian field which has the same operator self-similarity as in (1.4).

For the purpose of this paper, we will only be interested in the case $\alpha = H \in (0, 1)^N$. It follows from (1.2) that for every $t \in \mathbb{R}_+^N$,

$$B_0^H(t) = \kappa_H^{-1} X_0^H(t) + \kappa_H^{-1} \int_{(-\infty, t] \setminus [0, t]} \prod_{\ell=1}^N g_{H_\ell}(t_\ell, s_\ell) W(\mathrm{d}s), \tag{2.3}$$

and the two processes on the right-hand side of (2.3) are independent. We will show that in studying the regularity properties of the local times of $B^H$, the fractional Liouville sheet $X_0^H$ plays a crucial role and the second process in (2.3) can be neglected. More precisely, we will make use of the following property: For all integers $n \geq 2$, $t^1, \ldots, t^n \in \mathbb{R}_+^N$ and $u_1, \ldots, u_n \in \mathbb{R}$, we have

$$\mathrm{Var}\left(\sum_{j=1}^n u_j B_0^H(t^j)\right) \geq \kappa_H^{-2} \mathrm{Var}\left(\sum_{j=1}^n u_j X_0^H(t^j)\right). \tag{2.4}$$

Here and in the sequel, $\mathrm{Var}(\xi)$ denotes the variance of the random variable $\xi$.

Next we use an argument in [3] to provide a useful decomposition for $X_0^H(t)$. Let $\varepsilon > 0$ be fixed. For every $t \in [\varepsilon, \infty)^N$, we decompose the rectangle $[0, t]$ into the following disjoint union of sub-rectangles:

$$[0, t] = [0, \varepsilon]^N \cup \bigcup_{\ell=1}^N R_\ell(t) \cup \Delta(\varepsilon, t), \tag{2.5}$$

where $R_\ell(t) := R_\ell(\varepsilon, t) = \{r \in [0, t]^N \colon 0 \leq r_i \leq \varepsilon \text{ if } i \neq \ell, \varepsilon < r_\ell \leq t_\ell\}$ and $\Delta(\varepsilon, t)$ can be written as a union of $2^N - N - 1$ sub-rectangles of $[0, t]$. Denote the integrand in (2.2) by $g(t, r)$. It follows from (2.5) that for every $t \in [\varepsilon, \infty)^N$,

$$X_0^H(t) = \int_{[0,\varepsilon]^N} g(t, r) W(\mathrm{d}r)$$
$$+ \sum_{\ell=1}^N \int_{R_\ell(t)} g(t, r) W(\mathrm{d}r) + \int_{\Delta(\varepsilon, t)} g(t, r) W(\mathrm{d}r)$$
$$:= X(\varepsilon, t) + \sum_{\ell=1}^N Y_\ell(t) + Z(\varepsilon, t). \tag{2.6}$$

Since $\{X(\varepsilon, t), t \in [\varepsilon, \infty)^N\}$, $\{Y_\ell(t), t \in [\varepsilon, \infty)^N\}$ ($1 \leq \ell \leq N$) and $\{Z(\varepsilon, t), t \in [\varepsilon, \infty)^N\}$ are defined by the stochastic integrals w.r.t. $W$ over disjoint sets, they are independent Gaussian random fields.

The following lemma shows that every process $Y_\ell(t)$ has the property of strong local nondeterminism along the $\ell$th direction. It will be essential to our proofs.

**Lemma 2.1.** *Let $\ell \in \{1, 2, \ldots, N\}$ and let $I = [a, b] \in \mathcal{A}$ be a fixed interval. For any integer $n \geq 2$, $t^1, \ldots, t^n \in [a, b]$ such that*

$$t_\ell^1 \leq t_\ell^2 \leq \cdots \leq t_\ell^n,$$

*the following inequality for the conditional variance holds:*

$$\mathrm{Var}(Y_\ell(t^n) \mid Y_\ell(t^j) : 1 \leq j \leq n-1) \geq c_{2,2} |t_\ell^n - t_\ell^{n-1}|^{2H_\ell}, \tag{2.7}$$

*where $c_{2,2} > 0$ is a constant depending on $\varepsilon$, $I$ and $H$ only.*

**Proof.** Working in the Hilbert space setting, the conditional variance in (2.7) is the square of the $L^2(\mathbb{P})$-distance of $Y_\ell(t^n)$ from the subspace generated by $Y_\ell(t^j)$ ($1 \leq j \leq n-1$). Hence it is sufficient to show that there exists a constant $c_{2,2}$ such that for all $a_j \in \mathbb{R}$ ($j = 1, \ldots, n-1$),

$$\mathbb{E}\left( Y_\ell(t^n) - \sum_{j=1}^{n-1} a_j Y_\ell(t^j) \right)^2 \geq c_{2,2} |t_\ell^n - t_\ell^{n-1}|^{2H_\ell}. \tag{2.8}$$

However, by splitting $R_\ell(t^n)$ into two disjoint parts and using the independence, we derive that

$$\begin{aligned}
\mathbb{E}\left( Y_\ell(t^n) - \sum_{j=1}^{n-1} a_j Y_\ell(t^j) \right)^2 &\geq \mathbb{E}\left( \int_{R_\ell(t^n) \setminus R_\ell(t^{n-1})} g(t^n, r) W(\mathrm{d}r) \right)^2 \\
&= \int_0^\varepsilon \cdots \int_{t_\ell^{n-1}}^{t_\ell^n} \cdots \int_0^\varepsilon \prod_{k=1}^N (t_k^n - r_k)^{2H_k - 1} \,\mathrm{d}r \\
&\geq c_{2,2} |t_\ell^n - t_\ell^{n-1}|^{2H_\ell}.
\end{aligned} \tag{2.9}$$

This proves (2.8) and hence Lemma 2.1. □

The following lemma relates the fractional Brownian sheet $B_0^H$ to the independent Gaussian random fields $Y_\ell$ ($\ell = 1, \ldots, N$).

**Lemma 2.2.** *Let $I = [a, b] \in \mathcal{A}$ be a fixed interval. For every integer $n \geq 2$, $t^1, \ldots, t^n \in [a, b]$ and $u_1, \ldots, u_n \in \mathbb{R}$, we have*

$$\mathrm{Var}\left( \sum_{j=1}^n u_j B_0^H(t^j) \right) \geq \kappa_H^{-2} \sum_{\ell=1}^N \mathrm{Var}\left( \sum_{j=1}^n u_j Y_\ell(t^j) \right). \tag{2.10}$$

*Moreover, for every $k \in \{1, \ldots, N\}$ and positive numbers $p_1, \ldots, p_k \geq 1$ satisfying $\sum_{\ell=1}^k p_\ell^{-1} = 1$, we have*

$$\frac{1}{[\det\mathrm{Cov}(B_0^H(t^1), \ldots, B_0^H(t^n))]^{1/2}} \leq \prod_{\ell=1}^k \frac{c_{2,3}^n}{[\det\mathrm{Cov}(Y_\ell(t^1), \ldots, Y_\ell(t^n))]^{1/(2p_\ell)}}, \tag{2.11}$$

*where $\det\mathrm{Cov}(Z_1, \ldots, Z_n)$ denotes the determinant of the covariance matrix of the Gaussian random vector $(Z_1, \ldots, Z_n)$.*

**Proof.** The inequality (2.10) follows directly from (2.4), (2.6) and the independence of $Y_\ell$ ($\ell = 1, \ldots, N$). To prove (2.11), we note that for any positive definite $n \times n$ matrix $\Gamma$,

$$\int_{\mathbb{R}^n} \frac{[\det(\Gamma)]^{1/2}}{(2\pi)^{n/2}} \exp\left( -\frac{1}{2} x' \Gamma x \right) \mathrm{d}x = 1. \tag{2.12}$$

It follows from (2.10), (2.12) and the generalized Hölder inequality (see, e.g., [15], p. 140) that

$$\frac{1}{[\det\mathrm{Cov}(B_0^H(t^1), \ldots, B_0^H(t^n))]^{1/2}}$$

$$
\begin{aligned}
&= \frac{1}{(2\pi)^{n/2}} \int_{\mathbb{R}^n} \exp\left[-\frac{1}{2}\operatorname{Var}\left(\sum_{j=1}^n u_j B_0^H(t^j)\right)\right] du_1 \cdots du_n \\
&\leq \frac{1}{(2\pi)^{n/2}} \int_{\mathbb{R}^n} \exp\left[-c\sum_{\ell=1}^N \operatorname{Var}\left(\sum_{j=1}^n u_j Y_\ell(t^j)\right)\right] du_1 \cdots du_n \\
&\leq \frac{1}{(2\pi)^{n/2}} \prod_{\ell=1}^k \left\{\int_{\mathbb{R}^n} \exp\left[-c\operatorname{Var}\left(\sum_{j=1}^n u_j Y_\ell(t^j)\right)\right] du_1 \cdots du_n\right\}^{1/p_\ell} \\
&\leq \prod_{\ell=1}^k \frac{c_{2,4}^n}{[\det\operatorname{Cov}(Y_\ell(t^1),\ldots,Y_\ell(t^n))]^{1/(2p_\ell)}}.
\end{aligned}
\qquad(2.13)
$$

This yields (2.11) and the lemma is proved. □

## 3. Joint continuity of the local times

We start by briefly recalling some aspects of the theory of local times. For excellent surveys on local times of random and deterministic vector fields, we refer to [13] and [10].

Let $X(t)$ be a Borel vector field on $\mathbb{R}^N$ with values in $\mathbb{R}^d$. For any Borel set $T \subseteq \mathbb{R}^N$, the occupation measure of $X$ on $T$ is defined as the following measure on $\mathbb{R}^d$:

$$\mu_T(\bullet) = \lambda_N\{t \in T\colon X(t) \in \bullet\}.$$

If $\mu_T$ is absolutely continuous with respect to $\lambda_d$, we say that $X(t)$ has *local times* on $T$, and define its local times, $L(\bullet, T)$, as the Radon–Nikodým derivative of $\mu_T$ with respect to $\lambda_d$, i.e.,

$$L(x,T) = \frac{d\mu_T}{d\lambda_d}(x), \quad \forall x \in \mathbb{R}^d.$$

In the above, $x$ is the so-called *space variable*, and $T$ is the *time variable* of the local times. Sometimes, we write $L(x,t)$ in place of $L(x,[0,t])$. Note that if $X$ has local times on $T$ then for every Borel set $S \subseteq T$, $L(x, S)$ also exists.

By standard martingale and monotone class arguments, one can deduce that the local times have a version, still denoted by $L(x, T)$, such that it is a kernel in the following sense:

(i) For each fixed $S \in \mathcal{B}(T)$, where $\mathcal{B}(T)$ is the family of Borel subsets of $T$, the function $x \mapsto L(x,S)$ is Borel measurable in $x \in \mathbb{R}^d$.
(ii) For every $x \in \mathbb{R}^d$, $L(x,\cdot)$ is Borel measure on $\mathcal{B}(T)$.

Moreover, $L(x, T)$ satisfies the following *occupation density formula*: For every Borel set $T \subseteq \mathbb{R}^N$, and for every measurable function $f\colon \mathbb{R}^d \to \mathbb{R}_+$,

$$\int_T f(X(t))\,dt = \int_{\mathbb{R}^d} f(x) L(x,T)\,dx. \qquad(3.1)$$

See Theorems 6.3 and 6.4 in [13].

Suppose we fix a rectangle $I = \prod_{i=1}^N [a_i, a_i + h_i]$ in $\mathcal{A}$. Then, whenever we can choose a version of the local time, still denoted by $L(x, \prod_{i=1}^N [a_i, a_i + t_i])$, such that it is a continuous function of $(x, t_1, \ldots, t_N)$ $\in \mathbb{R}^d \times \prod_{i=1}^N [0, h_i]$, $X$ is said to have a *jointly continuous local time* on $I$. When a local time is jointly continuous, $L(x, \bullet)$ can be extended to be a finite Borel measure supported on the level set

$$X^{-1}(x) \cap I = \{t \in I\colon X(t) = x\}; \qquad(3.2)$$

see [1] for details. In other words, local times often act as a natural measure on the level sets of $X$. As such, they are useful in studying the various fractal properties of level sets and inverse images of the vector field $X$. In this regard, we refer to [6, 12, 27] and [30].

Berman [5, 6, 7] developed Fourier analytic methods for studying the existence and regularity of the local times of Gaussian processes. His methods were extended by Pitt [25] and Geman and Horowitz [13] to Gaussian random fields. Let $X = \{X(t), t \in \mathbb{R}^N\}$ be a Gaussian random field with values in $\mathbb{R}^d$. It follows from (25.5) and (25.7) in [13] (see also [14, 25]) that for all $x, y \in \mathbb{R}^d$, $T \in \mathcal{A}$ and all integers $n \geq 1$,

$$\mathbb{E}[L(x,T)^n] = (2\pi)^{-nd} \int_{T^n} \int_{\mathbb{R}^{nd}} \exp\left(-\mathrm{i} \sum_{j=1}^n \langle u^j, x \rangle \right)$$
$$\times \mathbb{E} \exp\left(\mathrm{i} \sum_{j=1}^n \langle u^j, X(t^j) \rangle \right) \mathrm{d}\overline{u}\,\mathrm{d}\overline{t} \tag{3.3}$$

and for all even integers $n \geq 2$,

$$\mathbb{E}[(L(x,T) - L(y,T))^n] = (2\pi)^{-nd} \int_{T^n} \int_{\mathbb{R}^{nd}} \prod_{j=1}^n [\mathrm{e}^{-\mathrm{i}\langle u^j, x \rangle} - \mathrm{e}^{-\mathrm{i}\langle u^j, y \rangle}]$$
$$\times \mathbb{E} \exp\left(\mathrm{i} \sum_{j=1}^n \langle u^j, X(t^j) \rangle \right) \mathrm{d}\overline{u}\,\mathrm{d}\overline{t}, \tag{3.4}$$

where $\overline{u} = (u^1, \ldots, u^n), \overline{t} = (t^1, \ldots, t^n)$, and each $u^j \in \mathbb{R}^d, t^j \in T \subset (0, \infty)^N$. In the coordinate notation we then write $u^j = (u_1^j, \ldots, u_d^j)$. These identities are also very useful for studying the local times of infinitely divisible random fields as well; see [10, 12] and [20].

Xiao and Zhang [33] have proved that if $d < \sum_{\ell=1}^N \frac{1}{H_\ell}$, then for all intervals $I \in \mathcal{A}$, $B^H$ has local times $\{L(x,I), x \in \mathbb{R}^d\}$ on $I$ and $L(\cdot, I) \in L^2(\mathbb{P} \times \lambda_d)$. In the following, we prove that under the same condition, the local time has a version that is jointly continuous in both space and time variables.

**Theorem 3.1.** *Let $B^H = \{B^H(t), t \in \mathbb{R}_+^N\}$ be a fractional Brownian sheet in $\mathbb{R}^d$ with index $H = (H_1, \ldots, H_N) \in (0,1)^N$. If $d < \sum_{\ell=1}^N \frac{1}{H_\ell}$, then for all intervals $I \in \mathcal{A}$, $B^H$ has a jointly continuous local time on $I$ almost surely.*

To prove Theorem 3.1 we will, similar to [12, 30, 33], first use the Fourier analytic arguments to derive estimates on the moments of the local times (see Lemmas 3.7 and 3.10) and then apply a multiparameter version of Kolmogorov continuity theorem (cf. [17]). The new ingredients in this paper are the "sectorial local nondeterministic" properties of fractional Brownian sheets proved in [29] and the results on fractional Liouville sheets proved in Section 2.

We will also make use of the following lemmas. Among them, Lemma 3.2 is proved in [29] and Lemma 3.3 is essentially due to Cuzick and DuPreez [9] (see also [19]).

**Lemma 3.2.** *Let $B_0^H = \{B_0^H(t), t \in \mathbb{R}_+^N\}$ be a fractional Brownian sheet in $\mathbb{R}$ with index $H = (H_1, \ldots, H_N) \in (0,1)^N$. Then for every $\varepsilon > 0$, there is a constant $c_{3,1} > 0$ such that for all integers $n \geq 2$, $t^1, \ldots, t^n \in [\varepsilon, \infty)^N$,*

$$\mathrm{Var}(B_0^H(t^n) \mid B_0^H(t^j), j \neq n) \geq c_{3,1} \sum_{\ell=1}^N \min\{|t_\ell^n - t_\ell^j|^{2H_\ell}, 0 \leq j \leq n-1\}, \tag{3.5}$$

*where $t_\ell^0 = 0$ for $\ell = 1, \ldots, N$.*

**Lemma 3.3.** *Let $Z_1, \ldots, Z_n$ be mean zero Gaussian variables which are linearly independent, then for any nonnegative Borel function $g: \mathbb{R} \to \mathbb{R}_+$,*

$$\int_{\mathbb{R}^n} g(v_1) \exp\left[-\frac{1}{2} \operatorname{Var}\left(\sum_{j=1}^n v_j Z_j\right)\right] dv_1 \cdots dv_n$$

$$= \frac{(2\pi)^{(n-1)/2}}{(\det \operatorname{Cov}(Z_1, \ldots, Z_n))^{1/2}} \int_{-\infty}^{\infty} g\left(\frac{v}{\sigma_1}\right) e^{-v^2/2} dv,$$

*where $\sigma_1^2 = \operatorname{Var}(Z_1 | Z_2, \ldots, Z_n)$ is the conditional variance of $Z_1$ given $Z_2, \ldots, Z_n$.*

The following technical lemma is essential in establishing the moment estimates for the local times $L(x,T)$. Since it may be of independent interest, we state it in a more general form than is needed in this paper.

**Lemma 3.4.** *For any $q \in [0, \sum_{\ell=1}^N H_\ell^{-1})$, let $\tau \in \{1, \ldots, N\}$ be the integer such that*

$$\sum_{\ell=1}^{\tau-1} \frac{1}{H_\ell} \leq q < \sum_{\ell=1}^{\tau} \frac{1}{H_\ell} \tag{3.6}$$

*with the convention that $\sum_{\ell=1}^0 \frac{1}{H_\ell} := 0$. Then there exists a positive constant $\delta_\tau \leq 1$ depending on $(H_1, \ldots, H_N)$ only such that for every $\delta \in (0, \delta_\tau)$, we can find $\tau$ real numbers $p_\ell \geq 1$ $(1 \leq \ell \leq \tau)$ satisfying the following properties:*

$$\sum_{\ell=1}^{\tau} \frac{1}{p_\ell} = 1, \quad \frac{H_\ell q}{p_\ell} < 1, \quad \forall \ell = 1, \ldots, \tau, \tag{3.7}$$

*and*

$$(1-\delta) \sum_{\ell=1}^{\tau} \frac{H_\ell q}{p_\ell} \leq H_\tau q + \tau - \sum_{\ell=1}^{\tau} \frac{H_\tau}{H_\ell}. \tag{3.8}$$

*Furthermore, if we denote $\alpha_\tau := \sum_{\ell=1}^{\tau} \frac{1}{H_\ell} - q > 0$, then for any positive number $\rho \in (0, \frac{\alpha_\tau}{2\tau})$, there exists an index $\ell_0 \in \{1, \ldots, \tau\}$ such that*

$$\frac{H_{\ell_0} q}{p_{\ell_0}} + 2H_{\ell_0} \rho < 1. \tag{3.9}$$

**Remark 3.5.** *It is important to note that the choice of the numbers $p_\ell \geq 1$ $(1 \leq \ell \leq \tau)$ depends on $\delta$. Moreover, it follows from the proof below that, except for the case of $\tau = 2$, we can always take $\delta_\tau = 1$.*

**Proof of Lemma 3.4.** First we prove (3.7) and (3.8). If (3.6) holds for $\tau = 1$, then for all $0 < \delta < \delta_1 := 1$, we can take $p_1 = 1$ and both (3.7) and (3.8) hold automatically.

We now prove the cases of $\tau \geq 2$ by induction. Our proof provides a general procedure for constructing a sequence $\{p_\ell, 1 \leq \ell \leq \tau\}$ of real numbers $p_\ell \geq 1$ satisfying (3.7) and (3.8) (there are many possible choices).

Assume that (3.6) holds for $\tau = 2$. We distinguish two cases: (i) $H_1 = H_2$ and (ii) $H_1 \neq H_2$. In the first case, we have $H_1^{-1} \leq q < 2H_1^{-1}$. We choose $\eta > 0$ such that

$$0 < \eta < \frac{(2 - H_1 q) H_1 q}{H_1 q - 1}$$

(if $H_1 q = 1$, then $\eta > 0$ can be arbitrarily chosen) and define

$$\frac{1}{p_1} = \frac{1}{H_1 q + \eta} \quad \text{and} \quad \frac{1}{p_2} = 1 - \frac{1}{p_1}.$$

Then a few lines of calculation verify that $p_1$ and $p_2$ satisfy (3.7) and (3.8) for all $\delta \in (0,1)$.

To consider the case (ii) we may and will assume, without loss of generality, that $H_1 < H_2$. Since $q < H_1^{-1} + H_2^{-1}$, there exists $\delta_2 > 0$ such that for all $\delta \in (0, \delta_2)$,

$$H_1 H_2 q(H_2 - H_1 + \delta H_1) < (H_2 - H_1)(H_2 + H_1 - \delta H_1). \tag{3.10}$$

For each fixed $\delta \in (0, \delta_2)$, we define

$$\frac{1}{p_1} = \frac{1}{1-\delta} \cdot \frac{1}{H_1 q} - \frac{\delta}{1-\delta} \cdot \frac{H_2}{H_2 - H_1} \quad \text{and} \quad \frac{1}{p_2} = 1 - \frac{1}{p_1}.$$

Then (3.7) follows from (3.6) and (3.10), and the equality sign in (3.8) holds.

Now we assume that the properties (3.7) and (3.8) hold for $\tau = n \in \{2, \ldots, N-1\}$ and consider the case of $\tau = n + 1$. Then we have

$$\sum_{\ell=2}^{n} \frac{1}{H_\ell} \leq q - \frac{1}{H_1} < \sum_{\ell=2}^{n+1} \frac{1}{H_\ell}. \tag{3.11}$$

Let $\delta \in (0,1)$ be fixed and we choose $\delta' \in (0, \delta \wedge \delta_n)$. Then it follows from (3.11) and the induction hypothesis that there exist $n$ constants $p'_\ell \geq 1$ ($\ell = 2, \ldots, n+1$) such that

$$\sum_{\ell=2}^{n+1} \frac{1}{p'_\ell} = 1, \qquad \frac{(q - 1/H_1) H_\ell}{p'_\ell} < 1, \quad \forall \ell = 2, \ldots, n+1, \tag{3.12}$$

and

$$(1 - \delta') \sum_{\ell=2}^{n+1} \frac{H_\ell(q - 1/H_1)}{p'_\ell} \leq H_{n+1}\left(q - \frac{1}{H_1}\right) + n - \sum_{\ell=2}^{n+1} \frac{H_{n+1}}{H_\ell}. \tag{3.13}$$

To define the constants $p_1, \ldots, p_{n+1}$ with the desired properties, we choose a constant $\eta > 0$ small so that

$$\frac{H_\ell q}{p'_\ell}\left(1 - \frac{1}{H_1 q} + \eta\right) < 1, \quad \forall \ell = 2, \ldots, n+1, \tag{3.14}$$

and

$$\frac{(1-\delta)(1 + (H_1 q\eta)/(H_1 q - 1))}{1 - \delta'} \leq 1. \tag{3.15}$$

This is possible because of (3.12).

Now we define $p_\ell$ ($1 \leq \ell \leq n+1$) by

$$\frac{1}{p_\ell} = \frac{1}{p'_\ell}\left(1 - \frac{1}{H_1 q} + \eta\right), \quad \forall \ell = 2, \ldots, n+1, \tag{3.16}$$

and

$$\frac{1}{p_1} = \frac{1}{H_1 q} - \eta. \tag{3.17}$$

It follows from this definition and (3.14) that

$$\sum_{\ell=1}^{n+1} \frac{1}{p_\ell} = 1 \quad \text{and} \quad \frac{H_\ell q}{p_\ell} < 1, \quad \forall \ell = 1, 2, \ldots, n+1. \tag{3.18}$$

That is, (3.7) holds for $\tau = n+1$. On the other hand, by some elementary calculation and (3.15) we can verify that

$$(1-\delta)\sum_{\ell=1}^{n+1}\frac{H_\ell q}{p_\ell} \leq H_{n+1}q + (n+1) - \sum_{\ell=1}^{n+1}\frac{H_{n+1}}{H_\ell}. \tag{3.19}$$

That is, (3.8) also holds for $\tau = n+1$. Hence the proof of (3.7) and (3.8) is completed.

Finally we prove (3.9). By (3.7), for every $\ell \in \{1,\ldots,\tau\}$, $\exists \varepsilon_\ell \in (0,1)$ such that $\frac{H_\ell q}{p_\ell} = 1 - \varepsilon_\ell$. Hence,

$$\sum_{\ell=1}^{\tau}\frac{\varepsilon_\ell}{H_\ell} = \sum_{\ell=1}^{\tau}\frac{1}{H_\ell} - \sum_{\ell=1}^{\tau}\frac{q}{p_\ell} = \sum_{\ell=1}^{\tau}\frac{1}{H_\ell} - q = \alpha_\tau > 0. \tag{3.20}$$

Hence there exists $\ell_0 \in \{1,\ldots,\tau\}$ such that $\varepsilon_{\ell_0} \geq \frac{H_{\ell_0}\alpha_\tau}{\tau}$. Note that for every positive number $\rho \in (0, \frac{\alpha_\tau}{2\tau})$, we have $2H_{\ell_0}\rho < \frac{H_{\ell_0}\alpha_\tau}{\tau} \leq \varepsilon_{\ell_0}$. Therefore

$$\frac{H_{\ell_0}q}{p_{\ell_0}} + 2H_{\ell_0}\rho = 1 - \varepsilon_{\ell_0} + 2H_{\ell_0}\rho < 1, \tag{3.21}$$

which completes the proof of (3.9). □

The following inequalities (3.22) and (3.23) with $a = 0$ are well known; see, e.g., [12]. The case $a > 0$ makes it possible for us to apply Lemma 3.4 for proving Lemmas 3.7 and 3.10.

**Lemma 3.6.** *For all integers $n \geq 1$, positive numbers $a$, $r$, $0 < b_j < 1$ and an arbitrary $s_0 \in [0, a/2]$,*

$$\int_{a \leq s_1 \leq \cdots \leq s_n \leq a+r} \prod_{j=1}^{n}(s_j - s_{j-1})^{-b_j}\,\mathrm{d}s_1\cdots\mathrm{d}s_n \leq c_{3,2}^n (n!)^{(1/n)\sum_{j=1}^{n}b_j - 1} r^{n - \sum_{j=2}^{n}b_j}, \tag{3.22}$$

*where $c_{3,2} > 0$ is a constant depending on $a$ and $b_j$'s only. In particular, if $b_j = \alpha$ for all $j = 1,\ldots,n$, then*

$$\int_{a \leq s_1 \leq \cdots \leq s_n \leq a+r} \prod_{j=1}^{n}(s_j - s_{j-1})^{-\alpha}\,\mathrm{d}s_1\cdots\mathrm{d}s_n \leq c_{3,2}^n (n!)^{\alpha-1} r^{n(1-(1-1/n)\alpha)}. \tag{3.23}$$

**Proof.** For simplicity, we only give the proof of (3.23) here. The proof of (3.22) is almost identical, and thus omitted. By integrating the integral in (3.23) in the order of $\mathrm{d}s_n, \mathrm{d}s_{n-1}, \ldots, \mathrm{d}s_1$, and by using a change of variable in each step to construct Beta functions, we derive

$$\int_{a \leq s_1 \leq \cdots \leq s_n \leq a+r} \prod_{j=1}^{n}(s_j - s_{j-1})^{-\alpha}\,\mathrm{d}s_1\cdots\mathrm{d}s_n$$
$$= \frac{1}{1-\alpha}\cdot\frac{\Gamma(2-\alpha)[\Gamma(1-\alpha)]^{n-2}}{\Gamma(1+(n-1)(1-\alpha))}\int_a^{a+r}(a+r-s_1)^{(n-1)(1-\alpha)}(s_1-s_0)^{-\alpha}\,\mathrm{d}s_1. \tag{3.24}$$

The inequality (3.23) follows from (3.24) and the Stirling's formula. □

In the rest of this section, we assume that $d < \sum_{\ell=1}^{N}\frac{1}{H_\ell}$ and $I \in \mathcal{A}$ is a fixed interval. For convenience, we further assume in the rest of this paper that

$$0 < H_1 \leq \cdots \leq H_N < 1. \tag{3.25}$$

We proceed to establish the moment estimates for the local times $L(x, T)$ which will be useful for proving the joint continuity of local times.

**Lemma 3.7.** *Let $B^H = \{B^H(t), t \in \mathbb{R}_+^N\}$ be a fractional Brownian sheet in $\mathbb{R}^d$ with index $H = (H_1, \ldots, H_N)$. If for some integer $\tau \in \{1, \ldots, N\}$ we have*

$$\sum_{\ell=1}^{\tau-1} \frac{1}{H_\ell} \leq d < \sum_{\ell=1}^{\tau} \frac{1}{H_\ell}, \tag{3.26}$$

*then there exists a positive and finite constant $c_{3,3}$, depending on $N, d, H$ and $I$ only, such that for all intervals $T = [a, a + \langle r \rangle] \subseteq I$ with edge-length $r \in (0, 1)$, all $x \in \mathbb{R}^d$ and all integers $n \geq 1$,*

$$\mathbb{E}[L(x,T)^n] \leq c_{3,3}^n (n!)^{N-\beta_\tau} r^{n\beta_\tau}, \tag{3.27}$$

*where $\beta_\tau = N - \tau - H_\tau d + \sum_{\ell=1}^{\tau} H_\tau/H_\ell$.*

**Remark 3.8.** *As we mentioned earlier, the local time $L(x, \bullet)$ may be extended as a random Borel measure supported on the level set $\Gamma_x = \{t \in (0, \infty)^N : B^H(t) = x\}$. Hence the moment estimate (3.27) contains a lot of information about the fractal properties of $\Gamma_x$. By Theorem 5 of [3], the Hausdorff dimension of the level set is given by*

$$\dim_H \Gamma_x = \min\left\{N - k - H_k d + \sum_{\ell=1}^{k} \frac{H_k}{H_\ell}, 1 \leq k \leq N\right\}, \tag{3.28}$$

*and the minimum is achieved by $\beta_\tau = N - \tau - H_\tau d + \sum_{\ell=1}^{\tau} H_\tau/H_\ell$, where $\tau$ satisfies (3.26). It is important to note that (3.27) is sharp and can be applied to strengthen the Hausdorff dimension result (3.28). We believe that the function $\varphi_1(r) = r^{\beta_\tau} (\log \log 1/r)^{N-\beta_\tau}$ is an exact Hausdorff measure function for $\Gamma_x$, and we will give a proof for the lower bound of the $\varphi_1$-Hausdorff measure of $\Gamma_x$ in Section 4. However, since the upper bound part relies on different methods, we will have to deal with it elsewhere.*

**Proof of Lemma 3.7.** For later use, we will start with an arbitrary closed interval $T = \prod_{\ell=1}^{N} [a_\ell, a_\ell + r_\ell] \subseteq I$. It follows from (3.3) and the fact that $B_1^H, \ldots, B_d^H$ are independent copies of $B_0^H$ that for all integers $n \geq 1$,

$$\mathbb{E}[L(x,T)^n] \leq (2\pi)^{-nd} \int_{T^n} \prod_{k=1}^{d} \left\{\int_{\mathbb{R}^n} \exp\left[-\frac{1}{2} \text{Var}\left(\sum_{j=1}^{n} u_k^j B_0^H(t^j)\right)\right] dU_k\right\} d\bar{t}, \tag{3.29}$$

where $U_k = (u_k^1, \ldots, u_k^n) \in \mathbb{R}^n$. Fix $k = 1, \ldots, d$ and denote the inner integral in (3.29) by $\mathcal{J}_k$. Then by Lemma 2.2, we have

$$\mathcal{J}_k \leq \int_{\mathbb{R}^n} \exp\left[-\frac{1}{2}\kappa_H^{-2} \sum_{\ell=1}^{N} \text{Var}\left(\sum_{j=1}^{n} u_k^j Y_\ell(t^j)\right)\right] dU_k$$

$$\leq \int_{\mathbb{R}^n} \exp\left[-\frac{1}{2}\kappa_H^{-2} \sum_{\ell=1}^{\tau} \text{Var}\left(\sum_{j=1}^{n} u_k^j Y_\ell(t^j)\right)\right] dU_k. \tag{3.30}$$

Since (3.26) holds, we apply Lemma 3.4 with $\delta = n^{-1}$ and $q = d$ to obtain $\tau$ positive numbers $p_1, \ldots, p_\tau \geq 1$ satisfying (3.7) and (3.8).

Applying the generalized Hölder inequality ([15], p. 140) to the last integral in (3.30), we derive that

$$\mathcal{J}_k \leq \prod_{\ell=1}^{\tau} \left\{\int_{\mathbb{R}^n} \exp\left[-\frac{p_\ell}{2}\kappa_H^{-2} \text{Var}\left(\sum_{j=1}^{n} u_k^j Y_\ell(t^j)\right)\right] dU_k\right\}^{1/p_\ell}$$

$$= c_{3,4}^n \prod_{\ell=1}^{\tau} [\det\text{Cov}(Y_\ell(t^1), \ldots, Y_\ell(t^n))]^{-1/(2p_\ell)}, \tag{3.31}$$

where the last equality follows from (2.12). Hence it follows from (3.29) and (3.31) that

$$\mathbb{E}[L(x,T)^n] \leq c_{3,5}^n \int_{T^n} \prod_{\ell=1}^{\tau} [\det\mathrm{Cov}(Y_\ell(t^1),\ldots,Y_\ell(t^n))]^{-d/(2p_\ell)} \, d\bar{t}. \tag{3.32}$$

To evaluate the integral in (3.32), we will first integrate $[dt_\ell^1 \cdots dt_\ell^n]$ for $\ell = 1,\ldots,\tau$. To this end, we will make use of the following fact about multivariate normal distributions: For any Gaussian random vector $(Z_1,\ldots,Z_n)$,

$$\det\mathrm{Cov}(Z_1,\ldots,Z_n) = \mathrm{Var}(Z_1) \prod_{j=2}^n \mathrm{Var}(Z_j|Z_1,\ldots,Z_{j-1}). \tag{3.33}$$

By the above fact and Lemma 2.1, we can derive that for every $\ell \in \{1,\ldots,\tau\}$ and for all $t^1,\ldots,t^n \in T = \prod_{\ell=1}^N [a_\ell, a_\ell + r_\ell]$ satisfying

$$a_\ell \leq t_\ell^{\pi_\ell(1)} \leq t_\ell^{\pi_\ell(2)} \leq \cdots \leq t_\ell^{\pi_\ell(n)} \leq a_\ell + r_\ell \tag{3.34}$$

for some permutation $\pi_\ell$ of $\{1,\ldots,N\}$, we have

$$\det\mathrm{Cov}(Y_\ell(t^1),\ldots,Y_\ell(t^n)) \geq c_{3,6}^n \prod_{j=1}^n (t_\ell^{\pi_\ell(j)} - t_\ell^{\pi_\ell(j-1)})^{2H_\ell}, \tag{3.35}$$

where $t_\ell^{\pi_\ell(0)} := \varepsilon$ (recall the decomposition (2.6)). We have chosen $\varepsilon < \frac{1}{2}\min\{a_\ell, 1 \leq \ell \leq N\}$ so that Lemma 3.6 is applicable.

It follows from (3.34) and (3.35) that

$$\int_{[a_\ell,a_\ell+r_\ell]^n} [\det\mathrm{Cov}(Y_\ell(t^1),\ldots,Y_\ell(t^n))]^{-d/(2p_\ell)} \, dt_\ell^1 \cdots dt_\ell^n$$

$$\leq \sum_{\pi_\ell} c^n \int_{a_\ell \leq t_\ell^{\pi_\ell(1)} \leq \cdots \leq t_\ell^{\pi_\ell(n)} \leq a_\ell + r_\ell} \prod_{j=1}^n \frac{1}{(t_\ell^{\pi_\ell(j)} - t_\ell^{\pi_\ell(j-1)})^{H_\ell d/p_\ell}} \, dt_\ell^1 \cdots dt_\ell^n$$

$$\leq c_{3,7}^n (n!)^{H_\ell d/p_\ell} r_\ell^{n(1-(1-1/n)H_\ell d/p_\ell)}. \tag{3.36}$$

In the above, the last inequality follows from (3.23).

Combining (3.32), (3.36) and continuing to integrate $[dt_\ell^1 \cdots dt_\ell^n]$ for $\ell = \tau + 1, \ldots, N$, we obtain

$$\mathbb{E}[L(x,T)^n] \leq c_{3,8}^n (n!)^{\sum_{\ell=1}^\tau H_\ell d/p_\ell} \prod_{\ell=1}^\tau r_\ell^{n(1-(1-1/n)H_\ell d/p_\ell)} \cdot \prod_{\ell=\tau+1}^N r_\ell^n. \tag{3.37}$$

Now we consider the special case when $T = [a, a + \langle r \rangle]$, i.e. $r_1 = \cdots = r_N = r$. Equations (3.37) and (3.8) with $\delta = n^{-1}$ and $q = d$ together yield

$$\mathbb{E}[L(x,T)^n] \leq c_{3,9}^n (n!)^{\sum_{\ell=1}^\tau H_\ell d/p_\ell} r^{n(N-(1-n^{-1})\sum_{\ell=1}^\tau H_\ell d/p_\ell)}$$

$$\leq c_{3,10}^n (n!)^{N-\beta_\tau} r^{n\beta_\tau}. \tag{3.38}$$

This proves (3.27). □

**Remark 3.9.** In the proof of Lemma 3.7, if we apply the generalized Hölder inequality to the first integral in (3.30) with $N$ positive numbers $p_1,\ldots,p_N$ defined by

$$p_\ell = \sum_{i=1}^N \frac{H_\ell}{H_i}, \quad \ell = 1,\ldots,N,$$

*then the above proof leading to (3.37) shows that the following inequality*

$$\mathbb{E}[L(x,T)^n] \leq c_{3,11}^n (n!)^{N\nu} \lambda_N(T)^{n(1-\nu)} \tag{3.39}$$

*holds for every interval $T \subset I$, where $\nu = d/(\sum_{\ell=1}^{N} H_\ell^{-1}) \in (0,1)$. We will apply this inequality in the proof of Theorem 3.1.*

Lemma 3.7 implies that for all $n \geq 1$, $L(x,T) \in L^n(\mathbb{R}^d)$ a.s. (see [13], p. 42). Our next lemma estimates the moments of the increments of $L(x,T)$ in $x$.

**Lemma 3.10.** *Assume (3.26) holds for some $\tau \in \{1,\ldots,N\}$. Then there exists a constant $c_{3,12}$, depending on $N, d, H$ and $I$ only, such that for all hypercubes $T = [a, a + \langle r \rangle] \subseteq I$, $x,y \in \mathbb{R}^d$ with $|x-y| \leq 1$, all even integers $n \geq 1$ and all $\gamma \in (0, 1 \wedge \frac{\alpha_\tau}{2\tau})$,*

$$\mathbb{E}[(L(x,T) - L(y,T))^n] \leq c_{3,12}^n (n!)^{N-\beta_\tau+(1+H_\tau)\gamma} |x-y|^{n\gamma} r^{n(\beta_\tau - H_\tau \gamma)}. \tag{3.40}$$

**Proof.** Let $\gamma \in (0, 1 \wedge \frac{\alpha_\tau}{2\tau})$ be a constant. Note that by the elementary inequalities

$$|e^{iu} - 1| \leq 2^{1-\gamma} |u|^\gamma \quad \text{for all } u \in \mathbb{R} \tag{3.41}$$

and $|u+v|^\gamma \leq |u|^\gamma + |v|^\gamma$, we see that for all $u^1, \ldots, u^n, x, y \in \mathbb{R}^d$,

$$\prod_{j=1}^{n} |e^{-i\langle u^j, x\rangle} - e^{-i\langle u^j, y\rangle}| \leq 2^{(1-\gamma)n} |x-y|^{n\gamma} {\sum}' \prod_{j=1}^{n} |u_{k_j}^j|^\gamma, \tag{3.42}$$

where the summation $\sum'$ is taken over all the sequences $(k_1, \ldots, k_n) \in \{1, \ldots, d\}^n$.

It follows from (3.4) and (3.42) that for every even integer $n \geq 2$,

$$\begin{aligned}
&\mathbb{E}[(L(x,T) - L(y,T))^n] \\
&\leq (2\pi)^{-nd} 2^{(1-\gamma)n} |x-y|^{n\gamma} \\
&\quad \times {\sum}' \int_{T^n} \int_{\mathbb{R}^{nd}} \prod_{m=1}^{n} |u_{k_m}^m|^\gamma \mathbb{E}\exp\left(-i\sum_{j=1}^{n} \langle u^j, B^H(t^j)\rangle\right) d\overline{u}\, d\overline{t} \\
&\leq c_{3,13}^n |x-y|^{n\gamma} {\sum}' \int_{T^n} d\overline{t} \\
&\quad \times \prod_{m=1}^{n} \left\{\int_{\mathbb{R}^{nd}} |u_{k_m}^m|^{n\gamma} \exp\left[-\frac{1}{2}\operatorname{Var}\left(\sum_{j=1}^{n} \langle u^j, B^H(t^j)\rangle\right)\right] d\overline{u}\right\}^{1/n},
\end{aligned} \tag{3.43}$$

where the last inequality follows from the generalized Hölder inequality.

Now we fix a vector $\overline{k} = (k_1, k_2, \ldots, k_n) \in \{1, \ldots, d\}^n$ and $n$ points $t^1, \ldots, t^n \in T$ such that $t_\ell^1, \ldots, t_\ell^n$ are all distinct for every $1 \leq \ell \leq N$ (the set of such points has full $(nN)$-dimensional Lebesgue measure). Let $\mathcal{M} = \mathcal{M}(\overline{k}, \overline{t}, \gamma)$ be defined by

$$\mathcal{M} = \prod_{m=1}^{n} \left\{\int_{\mathbb{R}^{nd}} |u_{k_m}^m|^{n\gamma} \exp\left[-\frac{1}{2}\operatorname{Var}\left(\sum_{j=1}^{n} \langle u^j, B^H(t^j)\rangle\right)\right] d\overline{u}\right\}^{1/n}. \tag{3.44}$$

Note that $B_\ell^H$ $(1 \leq \ell \leq N)$ are independent copies of $B_0^H$. By Lemma 3.2, the random variables $B_\ell^H(t^j)$ $(1 \leq \ell \leq N, 1 \leq j \leq n)$ are linearly independent. Hence Lemma 3.3 gives

$$\int_{\mathbb{R}^{nd}} |u_{k_m}^m|^{n\gamma} \exp\left[-\frac{1}{2}\operatorname{Var}\left(\sum_{j=1}^{n} \langle u^j, B^H(t^j)\rangle\right)\right] d\overline{u}$$

$$= \frac{(2\pi)^{(nd-1)/2}}{[\det \text{Cov}(B_0^H(t^1), \ldots, B_0^H(t^n))]^{d/2}} \int_{\mathbb{R}} \left(\frac{v}{\sigma_m}\right)^{n\gamma} e^{-v^2/2} \, dv$$

$$\leq \frac{c_{3,14}^n (n!)^\gamma}{[\det \text{Cov}(B_0^H(t^1), \ldots, B_0^H(t^n))]^{d/2}} \frac{1}{\sigma_m^{n\gamma}}, \tag{3.45}$$

where $\sigma_m^2$ is the conditional variance of $B_{k_m}^H(t^m)$ given $B_i^H(t^j)$ ($i \neq k_m$ or $i = k_m$ but $j \neq m$), and the last inequality follows from Stirling's formula.

Combining (3.44) and (3.45) we obtain

$$\mathcal{M} \leq \frac{c_{3,15}^n (n!)^\gamma}{[\det \text{Cov}(B_0^H(t^1), \ldots, B_0^H(t^n))]^{d/2}} \prod_{m=1}^{n} \frac{1}{\sigma_m^\gamma}. \tag{3.46}$$

For $\delta = 1/n$ and $q = d$, let $p_\ell$ ($\ell = 1, \ldots, \tau$) be the constants as in Lemma 3.4. Observe that, since $\gamma \in (0, \frac{\alpha_\tau}{2\tau})$, there exists an $\ell_0 \in \{1, \ldots, \tau\}$ such that

$$\frac{H_{\ell_0} d}{p_{\ell_0}} + 2H_{\ell_0} \gamma < 1. \tag{3.47}$$

It follows from (3.46) and Lemma 2.2 that

$$\mathcal{M} \leq c_{3,16}^n (n!)^\gamma \prod_{\ell=1}^{\tau} \frac{1}{[\det \text{Cov}(Y_\ell(t^1), \ldots, Y_\ell(t^n))]^{d/(2p_\ell)}} \prod_{m=1}^{n} \frac{1}{\sigma_m^\gamma}. \tag{3.48}$$

The second product in (3.48) will be treated as a "perturbation" factor and will be shown to be small when integrated. For this purpose, we use again the independence of the coordinate processes of $B^H$ and Lemma 3.2 to derive

$$\sigma_m^2 = \text{Var}(B_{k_m}^H(t^m) | B_{k_m}^H(t^j), j \neq m)$$

$$\geq c_{3,17}^2 \sum_{\ell=1}^{N} \min\{|t_\ell^m - t_\ell^j|^{2H_\ell} : j \neq m\}. \tag{3.49}$$

For any $n$ points $t^1, \ldots, t^n \in T$, let $\pi_1, \ldots, \pi_N$ be $N$ permutations of $\{1, 2, \ldots, n\}$ such that for every $1 \leq \ell \leq N$,

$$t_\ell^{\pi_\ell(1)} \leq t_\ell^{\pi_\ell(2)} \leq \cdots \leq t_\ell^{\pi_\ell(n)}. \tag{3.50}$$

Then, by (3.49) and (3.50) we have

$$\prod_{m=1}^{n} \frac{1}{\sigma_m^\gamma} \leq \prod_{m=1}^{n} \frac{1}{c_{3,18} \sum_{\ell=1}^{N} [(t_\ell^{\pi_\ell(m)} - t_\ell^{\pi_\ell(m-1)}) \wedge (t_\ell^{\pi_\ell(m+1)} - t_\ell^{\pi_\ell(m)})]^{H_\ell \gamma}}$$

$$\leq \prod_{m=1}^{n} \frac{1}{c_{3,18} [(t_{\ell_0}^{\pi_{\ell_0}(m)} - t_{\ell_0}^{\pi_{\ell_0}(m-1)}) \wedge (t_{\ell_0}^{\pi_{\ell_0}(m+1)} - t_{\ell_0}^{\pi_{\ell_0}(m)})]^{H_{\ell_0} \gamma}}$$

$$\leq c_{3,18}^{-n} \prod_{m=1}^{n} \frac{1}{(t_{\ell_0}^{\pi_{\ell_0}(m)} - t_{\ell_0}^{\pi_{\ell_0}(m-1)})^{q_{\ell_0}^m H_{\ell_0} \gamma}}, \tag{3.51}$$

for some $(q_{\ell_0}^1, \ldots, q_{\ell_0}^n) \in \{0, 1, 2\}^n$ satisfying $\sum_{m=1}^{n} q_{\ell_0}^m = n$ and $q_{\ell_0}^1 = 0$. That is, we will only need to consider the contribution of $\sigma_m$ in the $\ell_0$th direction.

So far we have obtained all the ingredients for bounding the integral in (3.43) and the rest of the proof is quite similar to the proof of Lemma 3.7. It follows from (3.48) and (3.51) that

$$\int_{T^n} \mathcal{M}(\overline{k},\overline{t},\gamma)\,d\overline{t} \leq c_{3,19}^n (n!)^\gamma \int_{T^n} \prod_{\ell=1}^\tau \frac{1}{[\det\mathrm{Cov}(Y_\ell(t^1),\ldots,Y_\ell(t^n))]^{d/(2p_\ell)}}$$
$$\times \prod_{m=1}^n \frac{1}{(t_{\ell_0}^{\pi_{\ell_0}(m)} - t_{\ell_0}^{\pi_{\ell_0}(m-1)})^{q_{\ell_0}^m H_{\ell_0}\gamma}}\,d\overline{t}. \tag{3.52}$$

To evaluate the above integral, we will first integrate $[dt_\ell^1 \cdots dt_\ell^n]$ for every $\ell = 1,\ldots,\tau$. Let us first consider $\ell = \ell_0$. By using Lemma 2.1, (3.33), (3.22) and, thanks to (3.47) and the nature of $q_{\ell_0}^m$, we see that

$$\int_{[a_{\ell_0}, a_{\ell_0}+r_{\ell_0}]^n} \frac{1}{[\det\mathrm{Cov}(Y_{\ell_0}(t^1),\ldots,Y_{\ell_0}(t^n))]^{d/(2p_{\ell_0})}} \times \prod_{m=1}^n \frac{1}{(t_{\ell_0}^{\pi_{\ell_0}(m)} - t_{\ell_0}^{\pi_{\ell_0}(m-1)})^{q_{\ell_0}^m H_{\ell_0}\gamma}}\,dt_{\ell_0}^1 \cdots dt_{\ell_0}^n \tag{3.53}$$

$$\leq \sum_{\pi_{\ell_0}} c_{3,20}^n \int_{a_{\ell_0} \leq t_{\ell_0}^{\pi_{\ell_0}(1)} \leq \cdots \leq t_{\ell_0}^{\pi_{\ell_0}(n)} \leq a_{\ell_0}+r_{\ell_0}}$$
$$\times \prod_{m=1}^n (t_{\ell_0}^{\pi_{\ell_0}(m)} - t_{\ell_0}^{\pi_{\ell_0}(m-1)})^{-(H_{\ell_0}d/p_{\ell_0} + q_{\ell_0}^m H_{\ell_0}\gamma)}\,dt_{\ell_0}^1 \cdots dt_{\ell_0}^n$$

$$\leq c_{3,21}^n (n!)^{H_{\ell_0}d/p_{\ell_0} + H_{\ell_0}\gamma} r_{\ell_0}^{n[1-(1-1/n)H_{\ell_0}d/p_{\ell_0} - H_{\ell_0}\gamma]}. \tag{3.54}$$

In the above, $t_{\ell_0}^{\pi_{\ell_0}(0)} = \varepsilon$ as in the proof of Lemma 3.7 and the last inequality follows from (3.22).

Meanwhile, recall that, for every $\ell \neq \ell_0$ ($\ell \in \{1,\ldots,\tau\}$), we have shown in (3.36) that

$$\int_{[a_\ell, a_\ell+r_\ell]^n} [\det\mathrm{Cov}(Y_\ell(t^1),\ldots,Y_\ell(t^n))]^{-d/(2p_\ell)}\,dt_\ell^1 \cdots dt_\ell^n$$
$$\leq c_{3,7}^n (n!)^{H_\ell d/p_\ell} r_\ell^{n(1-(1-1/n)H_\ell d/p_\ell)}. \tag{3.55}$$

Finally, we proceed to integrate $[dt_\ell^1 \cdots dt_\ell^n]$ for $\ell = \tau+1,\ldots,N$. It follows from the above that

$$\int_{T^n} \mathcal{M}(\overline{k},\overline{t},\gamma)\,d\overline{t} \leq c_{3,22}^n (n!)^{\sum_{\ell=1}^\tau H_\ell d/p_\ell + H_{\ell_0}\gamma + \gamma}$$
$$\times r_{\ell_0}^{n[1-(1-1/n)H_{\ell_0}d/p_{\ell_0} - H_{\ell_0}\gamma]} \times \prod_{\ell \neq \ell_0}^\tau r_\ell^{n[1-(1-1/n)H_\ell d/p_\ell]} \prod_{\ell=\tau+1}^N r_\ell^n. \tag{3.56}$$

In particular, if $r_1 = \cdots = r_N = r \leq 1$, we combine (3.43) and (3.56) to obtain

$$\mathbb{E}[(L(x,T) - L(y,T))^n]$$
$$\leq c_{3,23}^n |x-y|^{n\gamma} (n!)^{\sum_{\ell=1}^\tau H_\ell d/p_\ell + H_{\ell_0}\gamma + \gamma} \cdot r^{n(N-(1-1/n)\sum_{\ell=1}^\tau H_\ell d/p_\ell - H_{\ell_0}\gamma)}$$
$$\leq c_{3,24}^n (n!)^{N-\beta_\tau + (1+H_\tau)\gamma} |x-y|^{n\gamma} r^{n(\beta_\tau - H_\tau\gamma)}. \tag{3.57}$$

The last inequality follows from the fact that $H_{\ell_0} \leq H_\tau$ and Lemma 3.4. This finishes the proof of Lemma 3.10. □

Now we are ready to prove Theorem 3.1. It is similar to the proof of Theorem 4.1 in [33] and we include it for the sake of completeness.

**Proof of Theorem 3.1.** Let $I \in \mathcal{A}$ be fixed. For simplicity, we will assume $I = [\eta, 1]^N$ for some $\eta > 0$, say, $\eta = 2\varepsilon$ (cf. (2.6)). It follows from Lemma 3.10 and the multiparameter version of Kolmogorov's continuity theorem (cf. [17]) that, for every $T \in \mathcal{A}$ such that $T \subset I$, $B^H$ has almost surely a local time $L(x, T)$ that is continuous for all $x \in \mathbb{R}^d$.

To prove the joint continuity, observe that for all $x, y \in \mathbb{R}^d$ and $s, t \in I$, we have

$$\mathbb{E}[(L(x,[\eta,s]) - L(y,[\eta,t]))^n] \leq 2^{n-1}\{\mathbb{E}[(L(x,[\eta,s]) - L(x,[\eta,t]))^n] \\ + \mathbb{E}[(L(x,[\eta,t]) - L(y,[\eta,t]))^n]\}. \tag{3.58}$$

Since the difference $L(x,[\eta,s]) - L(x,[\eta,t])$ can be written as a sum of a finite number (only depends on $N$) of terms of the form $L(x, T_j)$, where each $T_j \in \mathcal{A}$ is a closed subinterval of $I$ with at least one edge length $\leq |s - t|$, we can use Lemma 3.7 and Remark 3.9, to bound the first term in (3.58). On the other hand, the second term in (3.58) can be dealt with using Lemma 3.10 as above. Consequently, for some $\gamma \in (0, 1)$ small, the right-hand side of (3.58) is bounded by $c_{3,25}^n (|x - y| + |s - t|)^{n\gamma}$, where $n \geq 2$ is an arbitrary even integer. Therefore the joint continuity of the local times follows again from the multiparameter version of Kolmogorov's continuity theorem. This finishes the proof of Theorem 3.1. □

We end this section with the following two technical lemmas, which will be useful in the next section.

**Lemma 3.11.** *Under the conditions of Lemma 3.7, there exist positive and finite constants $c_{3,26}$ and $c_{3,27}$, depending on $N, d, H$ and $I$ only, such that the following hold:*

(i) *For all $a \in I$ and hypercubes $T = [a, a + \langle r \rangle] \subseteq I$ with edge length $r \in (0,1)$, $x \in \mathbb{R}^d$ and all integers $n \geq 1$,*

$$\mathbb{E}[L(x + B^H(a), T)^n] \leq c_{3,26}^n (n!)^{N-\beta_\tau} r^{n\beta_\tau}, \tag{3.59}$$

*where $\beta_\tau = N - \tau - H_\tau d + \sum_{\ell=1}^{\tau} H_\tau / H_\ell$.*

(ii) *For all $a \in I$ and hypercubes $T = [a, a + \langle r \rangle] \subseteq I$, $x, y \in \mathbb{R}^d$ with $|x - y| \leq 1$, all even integers $n \geq 1$ and all $\gamma \in (0, 1 \wedge \frac{\alpha_\tau}{2\tau})$,*

$$\mathbb{E}[(L(x + B^H(a), T) - L(y + B^H(a), T))^n] \leq c_{3,27}^n (n!)^{N-\beta_\tau + (1+H_\tau)\gamma} |x - y|^{n\gamma} r^{n(\beta_\tau - H_\tau \gamma)}. \tag{3.60}$$

**Proof.** For each fixed $a \in I$, we define the Gaussian random field $Y = \{Y(t), t \in \mathbb{R}_+^N\}$ with values in $\mathbb{R}^d$ by $Y(t) = B^H(t) - B^H(a)$. It follows from (3.1) that if $B^H$ has a local time $L(x, S)$ on any Borel set $S$, then $Y$ also has a local time $\tilde{L}(x, S)$ on $S$ and, moreover, $L(x + B^H(a), S) = \tilde{L}(x, S)$. With little modification, the proofs of Lemmas 3.7 and 3.10 go through for the Gaussian field $Y$. Hence we derive that both (3.59) and (3.60) hold. □

The following lemma is a consequence of Lemma 3.11 and Chebyshev's inequality. The proof is standard, hence omitted.

**Lemma 3.12.** *Under the conditions of Lemma 3.7, there exist positive constants $c_{3,28}, c_{3,29}, b_1$ and $b_2 > 0$ (depending on $N, d, I$ and $H$ only), such that for all $a \in I$, $T = [a, a + \langle r \rangle]$ with $r \in (0, 1)$, $x \in \mathbb{R}^d$ and $u > 1$ large enough, we have*

$$\mathbb{P}\{L(x + B^H(a), T) \geq c_{3,28} r^{\beta_\tau} u^{N-\beta_\tau}\} \leq \exp(-b_1 u) \tag{3.61}$$

*and for $x, y \in \mathbb{R}^d$ with $|x - y| \leq 1$ and $\gamma \in (0, 1 \wedge \frac{\alpha_\tau}{2\tau})$,*

$$\mathbb{P}\{|L(x + B^H(a), T) - L(y + B^H(a), T)| \\ \geq c_{3,29} |x - y|^\gamma r^{\beta_\tau - H_\tau \gamma} u^{N-\beta_\tau + (1+H_\tau)\gamma}\} \leq \exp(-b_2 u). \tag{3.62}$$

## 4. Hölder conditions for the local times

In this section we investigate the local and uniform asymptotic behavior of the local time $L(x,T)$ at $x$ and the maximum local time $L^*(T) = \max_{x \in \mathbb{R}^d} L(x,T)$ as $\mathrm{diam}(T) \to 0$. The results are then applied to study the sample path properties of $B^H$.

*4.1. Hölder conditions for $L(x, \bullet)$*

By applying Lemma 3.12 (more precisely, (3.61) with $a = 0$) and the Borel–Cantelli lemma, one can easily derive the following law of the iterated logarithm for the local time $L(x, \cdot)$: If (3.26) holds for some $\tau \in \{1, \ldots, N\}$, then there exists a positive constant $c_{4,1}$ such that for every $x \in \mathbb{R}^d$ and $t \in (0, \infty)^N$,

$$\limsup_{r \to 0} \frac{L(x, U(t,r))}{\varphi_1(r)} \leq c_{4,1}, \tag{4.1}$$

where $U(t,r)$ is the open ball centered at $t$ with radius $r$ and $\varphi_1(r) = r^{\beta_\tau}(\log\log(1/r))^{N-\beta_\tau}$. It would be interesting to prove the lower bound in (4.1). For such a result for the local times of a one-parameter fractional Brownian motion, see [23].

It follows from Fubini's theorem that, with probability one, (4.1) holds for $\lambda_N$-almost all $t \in (0, \infty)^N$. Now we prove a stronger version of this result, which is useful in determining the exact Hausdorff measure of the level set.

**Theorem 4.1.** *Assume that $d < \sum_{\ell=1}^{N} \frac{1}{H_\ell}$. Let $\tau \in \{1, \ldots, N\}$ be the integer such that (3.26) holds and let $I \in \mathcal{A}$ be a fixed interval. For any fixed $x \in \mathbb{R}^d$, let $L(x, \cdot)$ be the local time of $B^H(t)$ at $x$ which is a random measure supported on the level set $(B^H)^{-1}(x)$. Then there exists a positive and finite constant $c_{4,2}$ independent of $x$ such that with probability $1$,*

$$\limsup_{r \to 0} \frac{L(x, U(t,r))}{\varphi_1(r)} \leq c_{4,2} \tag{4.2}$$

*holds for $L(x, \cdot)$-almost all $t \in I$, where $\varphi_1(r) = r^{\beta_\tau}(\log\log(1/r))^{N-\beta_\tau}$.*

**Proof.** The method of our proof is similar to that of Proposition 4.1 in [30]. For every integer $k > 0$, we consider the random measure $L_k(x, \bullet)$ on the Borel subsets $C$ of $I$ defined by

$$L_k(x, C) = \int_C (2\pi k)^{d/2} \exp\left(-\frac{k|B^H(t) - x|^2}{2}\right) dt$$
$$= \int_C \int_{\mathbb{R}^d} \exp\left(-\frac{|\xi|^2}{2k} + \mathrm{i}\langle \xi, B^H(t) - x\rangle\right) d\xi \, dt. \tag{4.3}$$

Then, by the occupation density formula (3.1) and the continuity of the function $y \mapsto L(y, C)$, one can verify that almost surely $L_k(x, C) \to L(x, C)$ as $k \to \infty$ for every Borel set $C \subset I$.

For every integer $m \geq 1$, denote $f_m(t) = L(x, U(t, 2^{-m}))$. From the proof of Theorem 3.1 we can see that almost surely the functions $f_m(t)$ are continuous and bounded. Hence we have almost surely, for all integers $m, n \geq 1$,

$$\int_I [f_m(t)]^n L(x, dt) = \lim_{k \to \infty} \int_I [f_m(t)]^n L_k(x, dt). \tag{4.4}$$

It follows from (4.3), (4.4) and the proof of Proposition 3.1 of [25] that for every positive integer $n \geq 1$,

$$\mathbb{E} \int_I [f_m(t)]^n L(x, dt) = \left(\frac{1}{2\pi}\right)^{(n+1)d} \int_I \int_{U(t, 2^{-m})^n} \int_{\mathbb{R}^{(n+1)d}} \exp\left(-\mathrm{i} \sum_{j=1}^{n+1} \langle x, u^j\rangle\right)$$

$$\times \mathbb{E} \exp\left(\mathrm{i} \sum_{j=1}^{n+1} \langle u^j, B^H(s^j)\rangle\right) \mathrm{d}\overline{u}\,\mathrm{d}\overline{s}, \tag{4.5}$$

where $\overline{u} = (u^1, \ldots, u^{n+1}) \in \mathbb{R}^{(n+1)d}$ and $\overline{s} = (t, s^1, \ldots, s^n)$. Similar to the proof of (3.27) we have that the right-hand side of Eq. (4.5) is at most

$$c_{4,3}^n \int_I \int_{U(t,2^{-m})^n} \frac{\mathrm{d}\overline{s}}{[\det\mathrm{Cov}(B_0^H(t), B_0^H(s^1), \ldots, B_0^H(s^n))]^{d/2}} \leq c_{4,4}^n (n!)^{N-\beta_\tau} 2^{-mn\beta_\tau}, \tag{4.6}$$

where $c_{4,4}$ is a positive finite constant depending on $N, d, H$, and $I$ only.

Let $\gamma > 0$ be a constant whose value will be determined later. We consider the random set

$$I_m(\omega) = \{t \in I: f_m(t) \geq \gamma \varphi_1(2^{-m})\}.$$

Denote by $\mu_\omega$ the restriction of the random measure $L(x, \cdot)$ on $I$, that is, $\mu_\omega(E) = L(x, E \cap I)$ for every Borel set $E \subset \mathbb{R}_+^N$. Now we take $n = \lfloor \log m \rfloor$, where $\lfloor y \rfloor$ denotes the integer part of $y$. Then by applying (4.6) and Stirling's formula, we have

$$\begin{aligned}
\mathbb{E}\mu_\omega(I_m) &\leq \frac{\mathbb{E} \int_I [f_m(t)]^n L(x, \mathrm{d}t)}{[\gamma \varphi_1(2^{-m})]^n} \\
&\leq \frac{c_{4,4}^n (n!)^{N-\beta_\tau} 2^{-mn\beta_\tau}}{\gamma^n 2^{-mn\beta_\tau} (\log m)^{n(N-\beta_\tau)}} \leq m^{-2},
\end{aligned} \tag{4.7}$$

provided $\gamma > 0$ is chosen large enough, say, $\gamma \geq c_{4,4}\mathrm{e}^2 := c_{4,2}$. This implies that

$$\mathbb{E}\left(\sum_{m=1}^\infty \mu_\omega(I_m)\right) < \infty.$$

Therefore, with probability 1 for $\mu_\omega$ almost all $t \in I$, we have

$$\limsup_{m\to\infty} \frac{L(x, U(t, 2^{-m}))}{\varphi_1(2^{-m})} \leq c_{4,2}. \tag{4.8}$$

Finally, for any $r > 0$ small enough, there exists an integer $m$ such that $2^{-m} \leq r < 2^{-m+1}$ and (4.8) is applicable. Since $\varphi_1(r)$ is increasing near $r = 0$, (4.2) follows from (4.8) and a monotonicity argument. □

**Theorem 4.2.** *Assume that $\sum_{\ell=1}^N \frac{1}{H_\ell} > d$ and $I \in \mathcal{A}$. Then there exists a positive constant $c_{4,5}$ such that for every $x \in \mathbb{R}^d$,*

$$\varphi_1\text{-}m((B^H)^{-1}(x) \cap I) \geq c_{4,5} L(x, I) \quad a.s., \tag{4.9}$$

*where $\varphi_1$-m denotes the $\varphi_1$-Hausdorff measure.*

**Proof.** As in the proof of Theorem 4.1 in [30], (4.9) follows from Theorem 4.1 and the upper density theorem of [26]. We omit the details. □

*4.2. Hölder conditions for $L^*(\bullet)$*

The following theorem establishes sharp Hölder conditions for the maximum local times $L^*(T) = \sup_{x \in \mathbb{R}^d} L(x, T)$ of fractional Brownian sheets as $\mathrm{diam}(T) \to 0$. Similar results for Brownian motion and some other random fields have been obtained by several authors. See, for example, [12, 16, 20, 30].

**Theorem 4.3.** *Let $B^H = \{B^H(t), t \in \mathbb{R}_+^N\}$ be a fractional Brownian sheet in $\mathbb{R}^d$ with index $H = (H_1, \ldots, H_N)$. We assume that there exists $\tau \in \{1, \ldots, N\}$ such that $H_1 = \cdots = H_\tau$ and $H_1 d < \tau$. Then, there exist positive constants $c_{4,6}$ and $c_{4,7}$ such that for every $s \in I$,*

$$\limsup_{r \to 0} \frac{L^*([s - \langle r \rangle, s + \langle r \rangle])}{r^{N - H_1 d} (\log \log r^{-1})^{H_1 d}} \leq c_{4,6} \quad a.s. \tag{4.10}$$

*and*

$$\limsup_{r \to 0} \sup_{s \in I} \frac{L^*([s - \langle r \rangle, s + \langle r \rangle])}{r^{N - H_1 d} (\log r^{-1})^{H_1 d}} \leq c_{4,7} \quad a.s. \tag{4.11}$$

For proving Theorem 4.3, we will make use of the following lemma, which is a consequence of Lemma 2.1 in [28] and Lemma 8 in [3].

**Lemma 4.4.** *Let $B^H = \{B^H(t), t \in \mathbb{R}_+^N\}$ be a fractional Brownian sheet in $\mathbb{R}^d$ with index $H = (H_1, \ldots, H_N)$ and let $I \in \mathcal{A}$ be fixed. Then there exist positive constants $c_{4,8}$ and $c_{4,9}$ such that for all $s \in I$, $T = [s, s + \langle h \rangle]$ with $h \in (0, 1)$ and all $u > c_{4,8} h^{H_1}$, we have*

$$\mathbb{P}\Big\{\sup_{t \in T} |B^H(t) - B^H(s)| \geq u\Big\} \leq \exp\left(-\frac{u^2}{c_{4,9} h^{2H_1}}\right). \tag{4.12}$$

**Proof of Theorem 4.3.** As in [12, 20, 30], the proof of Theorem 4.3 is based on Lemma 3.12 and a chaining argument. Hence we will only sketch a proof of (4.10), indicating the necessary modifications.

Let $g(r) = r^{N - H_1 d} (\log \log r^{-1})^{H_1 d}$ for $r > 0$ small enough. In order to prove (4.10) it is sufficient to show that for every $s \in I$,

$$\limsup_{n \to \infty} \frac{L^*(C_n)}{g(2^{-n})} \leq c_{4,10} \quad a.s., \tag{4.13}$$

where $C_n = [s, s + \langle 2^{-n} \rangle]$ for $n \geq 1$.

We divide the proof of (4.13) into four steps.

(a) Pick $u = 2^{-nH_1} \sqrt{2c_{4,9} \log n}$ in Lemma 4.4, we have

$$\mathbb{P}\Big\{\sup_{t \in C_n} |B^H(t) - B^H(s)| \geq 2^{-nH_1} \sqrt{2c_{4,9} \log n}\Big\} \leq \exp(-2 \log n) = n^{-2}. \tag{4.14}$$

Hence the Borel–Cantelli lemma implies that a.s. $\exists n_1 = n_1(\omega)$ such that

$$\sup_{t \in C_n} |B^H(t) - B^H(s)| \leq 2^{-nH_1} \sqrt{2c_{4,9} \log n}, \quad \text{for all } n \geq n_1. \tag{4.15}$$

(b) Let $\theta_n = 2^{-nH_1} (\log \log 2^n)^{-(1+H_1)}$ for all $n \geq 1$, and define

$$G_n = \{x \in \mathbb{R}^d \colon |x| \leq 2^{-nH_1} \sqrt{2c_{4,9} \log n} \text{ with } x = \theta_n p \text{ for some } p \in \mathbb{Z}^d\}.$$

Then, at least when $n$ is large enough, the cardinality of $G_n$ satisfies

$$\sharp G_n \leq c_{4,11} (\log n)^{(2+H_1)d}. \tag{4.16}$$

It follows from (3.61) that for any constant $c > 0$ and integer $n$ large enough,

$$\mathbb{P}\Big\{\max_{x \in G_n} L(x + B^H(s), C_n) \geq c^{H_1 d} g(2^{-n})\Big\} \leq c_{4,12} (\log n)^{(2+H_1)d} n^{-cb_1}. \tag{4.17}$$

(Note that $\beta_\tau = N - H_1 d$ under the assumptions of Theorem 4.3.) By choosing $c = 2b_1^{-1}$ in (4.17) we see that the right-hand side of (4.17) is summable. Hence, the Borel–Cantelli lemma implies that almost surely $\exists n_2 = n_2(\omega)$ such that

$$\max_{x \in G_n} L(x + B^H(s), C_n) \leq (2b_1^{-1})^{H_1 d} g(2^{-n}), \quad \text{for all } n \geq n_2. \tag{4.18}$$

(c) Given integers $n, k \geq 1$ and $x \in G_n$, we define

$$F(n, k, x) = \left\{ y \in \mathbb{R}^d \colon y = x + \theta_n \sum_{j=1}^k \varepsilon_j 2^{-j}, \varepsilon_j \in \{0,1\}^d \text{ for } 1 \leq j \leq k \right\}.$$

A pair of points $y_1, y_2 \in F(n, k, x)$ is said to be linked if $y_2 - y_1 = \theta_n \varepsilon 2^{-k}$ for some $\varepsilon \in \{0,1\}^d$. We choose $\gamma > 0$ small such that (3.62) in Lemma 3.12 holds, and then choose $\delta > 0$ such that $\delta(H_1 d + (1 + H_1)\gamma) < \gamma$. Consider the event $B_n$ defined by

$$B_n = \bigcup_{x \in G_n} \bigcup_{k=1}^\infty \bigcup_{y_1, y_2} \{|L(y_1 + B^H(s), C_n) - L(y_2 + B^H(s), C_n)|$$

$$\geq 2^{-n(N - H_1 d - H_1 \gamma)} |y_1 - y_2|^\gamma (c 2^{\delta k} \log n)^{H_1 d + (1 + H_1)\gamma} \}, \tag{4.19}$$

where $\bigcup_{y_1, y_2}$ signifies the union over all the linked pairs $y_1, y_2$, and where $c > 0$ is a constant whose value will be chosen later.

Note that $H_\tau = H_1$, by (3.62) we derive that for $n$ large enough,

$$\mathbb{P}\{B_n\} \leq c_{4,13} (\log n)^{(2+H_1)d} \sum_{k=1}^\infty 2^{(d+1)k} \exp(-cb_2 2^{\delta k} \log n)$$

$$\leq c_{4,14} (\log n)^{(2+H_1)d} n^{-cb_2}. \tag{4.20}$$

In the above the last inequality follows from the fact

$$\sum_{k=1}^\infty 2^{(d+1)k} \exp(-x 2^{\delta k}) \leq e^{-x}, \quad \forall x > 0 \text{ large enough.}$$

Hence, by choosing $c = 2b_2^{-1}$ in (4.19), the Borel–Cantelli lemma implies that almost surely, $B_n$ occurs only finitely many times.

(d) Fix an integer $n$ together with some $y \in \mathbb{R}^d$ that satisfies $|y| \leq 2^{-nH_1} \sqrt{2c_{4,9} \log n}$, we can represent $y$ in the form $y = \lim_{k \to \infty} y_k$ with

$$y_k = x + \theta_n \sum_{j=1}^k \varepsilon_j 2^{-j}, \tag{4.21}$$

where $y_0 = x \in G_n$ and $\varepsilon_j \in \{0,1\}^d$ for $j = 1, \ldots, k$.

Since the local time $L$ is jointly continuous, by expansion (4.21) and the triangular inequality, we see that on the event $B_n^c$,

$$|L(y + B^H(s), C_n) - L(x + B^H(s), C_n)|$$

$$\leq \sum_{k=1}^\infty |L(y_k + B^H(s), C_n) - L(y_{k-1} + B^H(s), C_n)|$$

$$\leq \sum_{k=1}^{\infty} 2^{-n(N-H_1 d - H_1 \gamma)} |y_k - y_{k-1}|^\gamma (2b_2^{-1} 2^{\delta k} \log n)^{H_1 d + (1+H_1)\gamma}$$

$$\leq c_{4,15} g(2^{-n}). \tag{4.22}$$

We combine (4.18) and (4.22) to get that for $n$ large enough,

$$\sup_{|x| \leq 2^{-nH_1}\sqrt{2c_{4,9}\log n}} L(x + B^H(s), C_n) \leq c_{4,16} g(2^{-n}). \tag{4.23}$$

That is

$$\sup_{|x - B^H(s)| \leq 2^{-nH_1}\sqrt{2c_{4,9}\log n}} L(x, C_n) \leq c_{4,16} g(2^{-n}). \tag{4.24}$$

Since $L^*(C_n) = \sup\{L(x, C_n): x \in \overline{B^H(C_n)}\}$, (4.13) follows from (4.24). This proves Theorem 4.3. □

The Hölder conditions for the local times of fractional Brownian sheets are closely related to the irregularity of the sample path of $B^H(t)$. To end this paper, we apply Theorem 4.3 to derive results about the degree of oscillation of the sample paths of $B^H(t)$, which greatly improves Theorem 3 of [3].

**Theorem 4.5.** *Let $B^H = \{B^H(t), t \in \mathbb{R}_+^N\}$ be an $(N, d)$-fractional Brownian sheet and let $I \in \mathcal{A}$ be a fixed interval. Then there exists a constant $c_{4,17} > 0$ such that for every $s \in I$,*

$$\liminf_{r \to 0} \sup_{t \in U(s,r)} \frac{|B^H(t) - B^H(s)|}{r^{H_1}(\log\log r^{-1})^{-H_1}} \geq c_{4,17} \quad a.s. \tag{4.25}$$

*and*

$$\liminf_{r \to 0} \inf_{t \in I} \sup_{t \in U(s,r)} \frac{|B^H(t) - B^H(s)|}{r^{H_1}(\log r^{-1})^{-H_1}} \geq c_{4,17} \quad a.s. \tag{4.26}$$

*In particular, the sample function $B^H(t)$ is almost surely nowhere differentiable in $(0, \infty)^N$.*

**Proof.** It is sufficient to prove the results for $d = 1$. Note that $H_1 < 1$, Theorem 4.3 is always applicable for $d = 1$ with $\tau = 1$. For any interval $Q \in \mathcal{A}$, we have

$$\lambda_N(Q) = \int_{\overline{B_0^H(Q)}} L(x, Q) \, dx \leq L^*(Q) \times \sup_{u,v \in Q} |B_0^H(u) - B_0^H(v)|. \tag{4.27}$$

By taking $Q = U(s, r)$ we see that (4.25) follows immediately from (4.27) and (4.10). Similarly, (4.26) follows from (4.27) and (4.11). □